\documentclass[tbtags,reqno]{amsart}
\usepackage{epsfig,amsthm,amsopn,amsfonts,amssymb}

\newtheorem{theorem}{Theorem}
\newtheorem{lemma}[theorem]{Lemma}
\newtheorem{proposition}[theorem]{Proposition}
\newtheorem{example}[theorem]{Example}

\newtheorem{corollary}[theorem]{Corollary}

\newtheorem{definition}[theorem]{Definition}

\newtheorem{property}[theorem]{Property}

\theoremstyle{remark}
\newtheorem{remark}[theorem]{Remark}

\def \CC {{\mathcal C}}
\def \CP {{\mathcal P}}

\def \la {\lambda}

\def\om {\omega}
\def\kbnd {\mathfrak c^{-1}}
\def\core {\mathfrak c}

\def \core {{\mathfrak c}}
\def\la {\lambda}
\def\gg {\gamma}
\def\aa {\alpha}

\def \dit#1{\text{\it ``#1"}}

\edef\savecatcodeat{\the\catcode`@}
\catcode`\@=11

\def\tb@ifSpecChars#1#2{#1}
\def\tb@ifNoSpecChars#1#2{#2}

\def\tableau{%
  \bgroup
  \@ifstar{\let\Tif\tb@ifNoSpecChars\tb@tableauB}
          {\let\Tif\tb@ifSpecChars\tb@tableauB}}

\def\tb@tableauB{
  \@ifnextchar[{\tb@tableauC}{\tb@tableauC[]}}

\def\tb@tableauC[#1]{\hbox\bgroup%
    \let\\=\cr
    \def\bl{\global\let\tbcellF\tb@cellNF}%
    \def\tf{\global\let\tbcellF\tb@cellH}
    \dimen2=\ht\strutbox \advance\dimen2 by\dp\strutbox%
    \ifx\baselinestretch\undefined\relax%
    \else%
       \dimen0=100sp \dimen0=\baselinestretch\dimen0%
       \dimen2=100\dimen2 \divide\dimen2 by\dimen0%
    \fi%
    \let\tpos\tb@vcenter
    \tb@initYoung
    \tb@options#1\eoo
    \let\arrow\tb@arrow%
    \dimen0=\Tscale\dimen2%
    \dimen1=\dimen0 \advance\dimen1 by \tb@fframe%
    \lineskip=0pt\baselineskip=0pt
    \def\tb@nothing{}%
    \def\endcellno{$\rss\egroup\bss\egroup}
    \def\endcell{\endcellno\kern-\dimen0}
    \def\begincell{\vbox to\dimen0\bgroup\vss\hbox to\dimen0\bgroup\hss$}%
    \let\overlay\tb@overlay%
    \let\fl\tb@fl%
    \let\lss\hss\let\rss\hss\let\tss\vss\let\bss\vss
    \def\mkcell##1{
        \let\tbcellF\tb@cellD
        \def\tb@cellarg{##1}
        \ifx\tb@cellarg\tb@nothing\let\tb@cellarg\tb@cellE\fi%
        \begincell\tb@cellarg\endcellno
        \tbcellF}
    \let\savecellF\tbcellF
     \Tif{\catcode`,=4\catcode`|=\active}{}\tb@tableauD}%

\let\tb@savetableauD\tableauD
{
    \catcode`|=\active \catcode`*=\active \catcode`~=\active%
    \catcode`@=\active
\gdef\tableauD#1{%
  \Tif{
    \mathcode`|="8000 \mathcode`*="8000%
    \mathcode`~="8000 \mathcode`@="8000%
    \def@{\bullet}%
    \let|\cr
    \let*\tf
    \let~\sk
  }{}%
  \tpos{\tabskip=0pt\halign{&\mkcell{##}\cr#1\crcr}}%
  \global\let\tbcellF\savecellF
  \egroup
  \egroup}
}
\let\tb@tableauD\tableauD
\let\tableauD\tb@savetableauD
\let\tb@savetableauD\undefined

\def\tb@options#1{\ifx#1\eoo\relax\else\tb@option#1\expandafter\tb@options\fi}

\def\tb@option#1{%
  \if#1t\let\tpos\tb@vtop\fi
  \if#1c\let\tpos\tb@vcenter\fi
  \if#1b\let\tpos\vbox\fi
  \if#1F\tb@initFerrers\fi
  \if#1Y\tb@initYoung\fi
  \if#1s\tb@initSmall\fi
  \if#1m\tb@initMedium\fi
  \if#1l\tb@initLarge\fi
  \if#1p\tb@initPartition\fi
  \if#1a\tb@initArrow\fi
}

\def\tb@vcenter#1{\ifmmode\vcenter{#1}\else$\vcenter{#1}$\fi}

\def\tb@vtop#1{\hbox{\raise\ht\strutbox\hbox{\lower\dimen0\vtop{#1}}}}

\def\tb@initPartition{\def\Tscale{.3}}
\def\tb@initSmall{\def\Tscale{1}}
\def\tb@initMedium{\def\Tscale{2}}
\def\tb@initLarge{\def\Tscale{3}}

\def\tb@initArrow{\dimen2=1.25em}

\def\tb@initYoung{%
  \def\tb@cellE{}
  \let\tb@cellD\tb@cellN
  \def\sk{\global\let\tbcellF\tb@cellNF}}
\def\tb@initFerrers{%
  \def\tb@cellE{\bullet}
  \let\tb@cellD\tb@cellNF
  \def\sk{\bullet}}

\tb@initMedium

\def\tb@sframe#1{%
  \vbox to0pt{
    \vss
    \hbox to0pt{%
      \hss
      \vbox to\dimen1{
        \hrule depth #1 height0pt
        \vss
        \hbox to\dimen1{
          \vrule width #1 height\dimen1
          \hss
          \vrule width #1
          }%
        \vss
        \hrule height #1 depth 0in
        }%
      \kern-\tb@hframe
      }%
    \kern-\tb@hframe}}

\def\tb@hframe{.2pt}\def\tb@fframe{.4pt}\def\tb@bframe{2pt}
\def\tb@cellH{\tb@sframe{\tb@bframe}}       
\def\tb@cellNF{}                            
\def\tb@cellN{\tb@sframe{\tb@fframe}}       
\let\tbcellF\tb@cellN                       

\def\tb@rpad{1pt}
\def\tb@lpad{1pt}
\def\tb@tpad{1.8pt}
\def\tb@bpad{1.8pt}

\def\tb@overlay{\endcell\@ifnextchar[{\tb@overlaya}{\begincell}}
\def\tb@overlaya[#1]{\vbox to\dimen0\bgroup%
  \tb@overlayoptions#1\eoo%
  \tss\hbox to\dimen0\bgroup\lss$}
\def\tb@overlayoptions#1{\ifx#1\eoo\relax\else\tb@overlayoption#1\expandafter\tb@overlayoptions\fi}

\def\tb@overlayoption#1{
  \if#1t\def\tss{\vskip\tb@tpad}\let\bss\vss\fi
  \if#1c\let\tss\vss\let\bss\vss\fi
  \if#1b\def\bss{\vskip\tb@bpad}\let\tss\vss\fi
  \if#1l\def\lss{\hskip\tb@lpad}\let\rss\hss\fi
  \if#1m\let\lss\hss\let\rss\hss\fi
  \if#1r\def\rss{\hskip\tb@rpad}\let\lss\hss\fi
}

\def\tb@fl{\endcell\begincell\vrule depth 0pt width \dimen0 height \dimen0 \endcell\begincell}

\@ifundefined{diagram}{}{
\def\tb@arrowpad{.5}

\newoptcommand{\tb@arrow}{\@ne}[2]{%
  \endcell
   \begingroup%
   \let\dg@getnodesize\tb@getnodesize
   \dg@USERSIZE=#1\relax%
   \ifnum\dg@USERSIZE<\@ne \dg@USERSIZE=\@ne \fi%
   \dg@parse{#2}%
   \dg@label{\tb@draw{#1}{#2}}}

\def\tb@getnodesize#1#2#3#4#5{\dimen3=\tb@arrowpad\dimen2 #4=\dimen3 #5=\dimen3\relax}
\def\tb@getnodesize#1#2#3#4#5{\ifnum#2=0\ifnum#3=0\tb@getnodesizetail{#4}{#5}\else\tb@getnodesizehead{#4}{#5}\fi\else\tb@getnodesizehead{#4}{#5}\fi}
\def\tb@getnodesizetail#1#2{\dimen3=.5\dimen2 #1=\dimen3 #2=\dimen3}
\def\tb@getnodesizehead#1#2{\dimen3=.5\dimen2 #1=\dimen3 #2=\dimen3}

\def\tb@draw#1#2#3#4{%
        \dg@X=0\dg@Y=0\dg@XGRID=1\dg@YGRID=1\unitlength=.001\dimen0%
        \dg@LBLOFF=\dgLABELOFFSET \divide\dg@LBLOFF\unitlength%
        \dg@drawcalc
        \begincell
        \let\lams@arrow\tb@lams@arrow
        \begin{picture}(0,0)\begingroup\dg@draw{#1}{#2}{#3}{#4}\end{picture}%
        \endcell
        \endgroup
        \begincell}
}

\def\tb@lams@arrow#1#2{%
 \lams@firstx\z@\lams@firsty\z@
 \lams@lastx#1\relax\lams@lasty#2\relax
 \lams@center\z@
 \N@false\E@false\H@false\V@false
 \ifdim\lams@lastx>\z@\E@true\fi
 \ifdim\lams@lastx=\z@\V@true\fi
 \ifdim\lams@lasty>\z@\N@true\fi
 \ifdim\lams@lasty=\z@\H@true\fi
 \NESW@false
 \ifN@\ifE@\NESW@true\fi\else\ifE@\else\NESW@true\fi\fi
 \ifH@\else\ifV@\else
  \lams@slope
  \ifnum\lams@tani>\lams@tanii
   \lams@ht\ten@\p@\lams@wd\ten@\p@
   \multiply\lams@wd\lams@tanii\divide\lams@wd\lams@tani
  \else
   \lams@wd\ten@\p@\lams@ht\ten@\p@
   \divide\lams@ht\lams@tanii\multiply\lams@ht\lams@tani
  \fi
 \fi\fi
 \ifH@  \lams@harrow
 \else\ifV@ \lams@varrow
 \else \lams@darrow
 \fi\fi
}

\catcode`\@=\savecatcodeat
\let\savecatcodeat\undefined

\begin{document}

\title[A $k$-tableau characterization of $k$-Schur functions]
{A $k$-tableau characterization of $k$-Schur functions}

\author{Luc Lapointe}
\thanks{Research supported in part by FONDECYT (Chile) grant \#1030114,
the Programa Formas Cuadr\'aticas of the Universidad de Talca,
and NSERC (Canada) grant \#250904}
\address{Instituto de Matem\'atica y F\'{\i}sica,
Universidad de Talca, Casilla 747, Talca, Chile}
\email{lapointe@inst-mat.utalca.cl}

\author{Jennifer Morse}
\thanks{Research supported in part by NSF grant \#DMS-0400628}
\address{Department of Mathematics,
University of Miami, Coral Gables, Fl 33124}
\email{morsej@math.miami.edu}


\begin{abstract}
We study $k$-Schur functions characterized by $k$-tableaux, proving 
combinatorial properties such as a $k$-Pieri rule and a $k$-conjugation.  
This new approach relies on developing the theory of $k$-tableaux, and
includes the introduction of a weight-permuting involution on these
tableaux that generalizes the Bender-Knuth involution.  This work lays the
groundwork needed to prove that the set of $k$-Schur Littlewood-Richardson
coefficients contains the 3-point Gromov-Witten invariants; structure
constants for the quantum cohomology ring. 
\end{abstract}

\maketitle


\section{Introduction}
The Schur functions $s_\lambda$ form a basis for the symmetric 
function space $\Lambda$ which plays a fundamental role in combinatorics, 
representation theory and algebraic geometry.  For instance, the 
Pieri formula for multiplying Schubert varieties in the intersection
ring of a Grassmannian is equivalent to the formula for multiplying a 
Schur function and a  homogeneous function $h_\ell$ in $\Lambda$:
\begin{equation}
\label{pieri}
h_\ell \, s_\mu= \sum_{\lambda=\mu+horizontal~strip} s_\lambda
\,.
\end{equation}
A formula defined by vertical strips, rather than horizontal,
describes the product of an elementary symmetric function
$e_\ell s_\mu$.
More generally, structure constants of the cohomology ring of the 
Grassmannian in the basis of Schubert classes are none other than the 
\dit{Littlewood-Richardson coefficients}, occurring in the expansion
\begin{equation}
\label{litric}
s_\nu\, s_\mu = 
\sum_{\lambda}
c_{\nu\mu}^\lambda s_\lambda
\,.
\end{equation}

\medskip

Combinatorics is deeply intertwined with the theory of Schur functions.  
The Littlewood-Richardson coefficients are characterized by certain 
skew tableaux and at a more fundamental level, the very definition 
of column-strict tableaux arises by iterating \eqref{pieri}.  That is,
\begin{equation}
\label{kostka}
h_\mu= \sum_{\lambda \trianglerighteq \mu} K_{\lambda \mu}\, s_\lambda
\,,
\end{equation}
where the \dit{Kostka numbers} $K_{\lambda\mu}$ count the number of
tableaux of shape $\lambda$ and weight $\mu$, with the column-strict
condition required by the Pieri rule.  
The role of Schur functions also ties into the combinatorial theory 
of partitions as can be seen when working with the algebra 
endomorphism defined by $\omega e_\ell = h_\ell$.  This 
involution acts simply on a Schur function by
\begin{equation} \label{invo}
\omega s_\lambda = s_{\lambda'}\,,
\end{equation}
where $\lambda'$ is the partition conjugate to $\lambda$.

\medskip

Recent developments in symmetric function theory involved the study 
of Macdonald polynomials.  The Schur basis is again fundamental in 
this setting since the Macdonald expansion coefficients 
in this basis have a representation 
theoretic interpretation \cite{[Ga],[Ga2],[Ha]}.
In work with Lascoux on the Macdonald polynomials \cite{[LLM]}, we 
discovered a new family of symmetric functions defined for each
partition $\lambda$, where $\lambda_1\leq k$, by
$$
s_\lambda^{(k)}=\sum_{T\in\mathcal S_\lambda} s_{\text{shape}(T)}\,,
$$
for certain sets of tableaux $\mathcal S_\lambda$.  Experimentation 
suggested that these functions play the fundamental role of the Schur 
functions in the subring $\Lambda^{(k)}=\mathbb Z[h_1,\ldots,h_k]$; 
they form a basis for $\Lambda^{(k)}$ that satisfies
generalizations of classical Schur function properties such as 
\eqref{pieri}, \eqref{litric}, 
\eqref{kostka} and \eqref{invo}.  As such, we coined the functions
\dit{$k$-Schur functions}.
Unfortunately, while fertile for intuition and computer experimentation, 
the characterization of $\mathcal S_\lambda$ lagged in mechanisms of 
proof.  This led us to seek an alternative characterization and in
\cite{[LMfil]}, we introduced functions that are conjecturally equivalent.  
Although we were able to prove that these functions form a basis for 
$\Lambda^{(k)}$, the combinatorial conjectures 
remained open.

\medskip

Continued empirical study of the $k$-Schur functions led us to
tangential results in algebraic combinatorics.  In particular, 
a new family of tableaux were drawn from the conjectured Pieri rule 
for $k$-Schurs:
\begin{equation}
\label{intropieri}
h_\ell \,s_{\mu}^{(k)} = \sum_{\lambda\in H_{\mu,\ell}^k} s_\lambda^{(k)}
\,,
\end{equation}
where $H_{\mu,\ell}^k$ is a certain subset of the partitions obtained
by adding a horizontal $\ell$-strip to  $\mu$.
Iterating this relation gives
\begin{equation} 
\label{hexp}
h_\mu= \sum_{\lambda} K_{\lambda\mu}^{(k)}\,
s_\lambda^{(k)} \,.
\end{equation} 
Using the Pieri rule as a guide, we defined the \dit{$k$-tableaux}
(see Definition~\ref{defktabgen}) as certain fillings of $k+1$-cores
whose enumeration gives the \dit{$k$-Kostka numbers} 
$K_{\lambda\mu}^{(k)}$.  In \cite{[LMcore]}, we proved that these 
tableaux directly connect to the type-A affine Weyl group, and explained 
their role in enumerating the monomial terms of coefficients in 
the $k$-Schur expansion of Macdonald polynomials.  

\medskip

The family of $k$-tableaux is the central object of study here.
As with usual tableaux, these are associated to a shape $\lambda$ 
and weight $\mu$ and satisfy \cite{[LMcore]}:
\begin{equation} \label{triangu}
K_{\lambda\mu}^{(k)} = 0 {\rm~~~when~~~} \lambda \ntrianglerighteq \mu\,
\quad\text{and}\quad
K_{\lambda\lambda}^{(k)} = 1\,.
\end{equation}
Thus, \eqref{hexp} gives an invertible system that can be used
to characterize the $k$-Schur functions.
In this paper we investigate this third (conjecturally equivalent) 
characterization for the $k$-Schur functions.  The very definition 
implies that these functions  form a basis for $\Lambda^{(k)}$.  
Then, with an in-depth study of $k$-tableaux, we are able to prove 
that these polynomials satisfy several combinatorial properties 
including analogs of \eqref{pieri}, \eqref{kostka}, and \eqref{invo}.  
Moreover, our results strongly suggest that $k$-tableaux are
the objects to approach two long-standing open problems -- finding 
a combinatorial interpretation for the 3-point Gromov-Witten 
invariants and for the Macdonald expansion coefficients.

\medskip

To be more specific, by proving a number of results about the structure
of $k$-tableaux, we discover an involution on the weight that
reduces to the Bender-Knuth involution \cite{[BK]} on column-strict
tableaux.  Consequently, we can derive the following relation on
$k$-Kostka numbers:
$$
K_{\lambda\alpha}^{(k)} = K_{\lambda\mu}^{(k)}
$$
for $\alpha$ any rearrangement of $\mu$.  From this, we
are able to prove the $k$-Pieri rule for $k$-Schur functions
\eqref{intropieri}.  We also prove a formula for $e_\ell s_\mu^{(k)}$
that depends on a subset of the shapes obtained by adding vertical 
$\ell$-strips to $\mu$.  From the $e_\ell$ and $h_\ell$ $k$-Pieri 
rules we can show that applying the $\omega$-involution to $s_\lambda^{(k)}$ 
produces exactly one $k$-Schur function:
\begin{equation}
\omega s_{\lambda}^{(k)} = s_{\lambda^{\omega_k}}^{(k)} \, ,
\end{equation}
indexed by the \dit{$k$-conjugate} of $\lambda$.  We show that 
$s_{\lambda}^{(k)}$ coincides with $s_{\lambda}$ when the hook-length
of $\lambda$ is not larger than $k$, and thus that $k$-Schur functions
reduce to Schur functions when $k$ is large enough.
This concurs with our assertion that the $k$-Schur functions 
are the \dit{Schur basis} for $\Lambda^{(k)}$, since $\Lambda^{(k)}=\Lambda$ 
when $k\to\infty$.  

\medskip

As mentioned, the $k$-tableaux are connected to the affine symmetric
group $\tilde S_{k+1}$.  In particular, the $k$-Pieri rule induces
an order that is isomorphic to the weak order on $\tilde S_{k+1}$
modulo a maximal parabolic subgroup isomorphic to $S_{k+1}$.  
The interpretation of the weak order on $\tilde S_{k+1}/S_{k+1}$ 
as the tiling of a cone in $k$-space by permutahedra can be 
seen on the level of symmetric functions by identifying 
vertices with $k$-Schur functions.  The vectors 
of translation invariance in the tiling turn out simply to be
usual Schur functions indexed by \dit{$k$-rectangles} -- partitions of 
the form $(\ell^{k-\ell+1})$.  This follows from our last property:
\begin{equation}
s_{\square} \, s_{\lambda}^{(k)} = s_{\lambda \cup \square }^{(k)} \, ,
\end{equation}
where $\square$ is any $k$-rectangle.  This result implies that there 
are $k!$ \dit{$k$-irreducible} $k$-Schur functions from which any 
other $k$-Schur can be constructed by multiplication with Schur functions 
indexed by $k$-rectangles.  These $k$-irreducibles are indexed by partitions 
with no more than $i$ parts equal to $k-i$.  This property, as well as its
$t$-generalization, also holds for the functions introduced in 
\cite{[LMfil]} (see \cite{[LMadv]}).

\medskip

Although this article concentrates on proving that the $k$-Schur 
functions are the fundamental combinatorial analog for the Schur 
functions in the subspace $\Lambda^{(k)}$,  this analogy extends
beyond combinatorics.  Results presented here are the tools needed 
to carry out the first step in this direction.
In \cite{[LMhecke]}, we prove that the $k$-Schur functions 
provide the natural basis for the quantum cohomology of the 
Grassmannian \cite{[Ag],[Wi]}.  Consequently, the three point 
Gromov-Witten invariants are none other than 
\dit{$k$-Littlewood-Richardson coefficients} occurring in
\begin{equation}
s_{\mu}^{(k)} \, s_{\nu}^{(k)} = \sum_{\lambda} \, c_{\mu \nu}^{\lambda, k} 
s_{\lambda}^{(k)} \, .
\end{equation}
This implies the positivity of $c_{\lambda\mu}^{\nu,k}$ in certain
cases; conjectured to hold in general. Explicit connections are 
also made in \cite{[LMhecke]} between $k$-Littlewood Richardson 
coefficients, fusion coefficients for the WZW-conformal field theories
\cite{[TUY]}, and structure constants related to certain representations 
of Hecke algebras at roots of unity \cite{[GW]}.

\medskip

Given these developments on $k$-tableaux and $k$-Schur functions, 
there are many natural paths for future work.  Most notable is to 
investigate a likely connection between the $k$-Schur functions and the 
affine (loop) Grassmannian.  In particular,  M. Shimozono conjectured 
that the $k$-Littlewood-Richardson coefficients give the integral homology 
of the loop Grassmannian.  There is extensive computational evidence that the 
\dit{dual $k$-Schur functions}, defined in \cite{[LMhecke]} by summing 
over the monomial weights of $k$-tableaux, are the Schubert classes
in the cohomology of the loop Grassmannian.  
Another topic to be explored is the problem of finding appropriate 
skew $k$-tableaux to combinatorially describe the $k$-Littlewood-Richardson 
coefficients (and consequently the 3-point Gromov-Witten invariants).
A last example goes back to the Macdonald problem from where the
$k$-Schur functions arose.  \eqref{kostka} gives in particular that
\begin{equation}
h_{1^n}\, =\, \sum_{\lambda} K_{\lambda 1^n}^{(k)}\, s_\lambda^{(k)}\,,
\end{equation}
where $K_{\lambda 1^n}$ enumerates the subset of \dit{standard} $k$-tableaux.
Thus, there should exist a pair of statistics on $k$-standard tableaux 
to explain the $k$-Schur expansion coefficients in a Macdonald
polynomial $H_\mu[X;q,t]$  since $h_{1^n}=H_\mu[X;1,1]$.  Exciting 
mathematics has sprung from the search for combinatorial interpretations
of the Gromov-Witten invariants and the Macdonald coefficients, 
however the conjectures remain open.  See 
\cite{[BMW],[BKMW],[BKT],[KMSW],[KTW],[SS1],[SS2],[Tu]} 
and \cite{[Hag],[HH]} for examples of recent progress in these directions.

\section{Definitions}
Let $\Lambda$ denote the ring of symmetric functions, generated by the 
elementary symmetric functions $e_r=\sum_{i_1<\ldots <i_r}x_{i_1}
\cdots x_{i_r}$, or equivalently by the complete symmetric functions
$h_r=\sum_{i_1\leq\ldots\leq i_r}x_{i_1}\cdots x_{i_r}$, and let
$\Lambda^{k} = \mathbb Z[h_1,\ldots,h_k]$.  Bases for $\Lambda$ are 
indexed by partitions $\lambda=(\lambda_1\geq\dots\geq\lambda_m>0)$ 
whose degree $\lambda$ is $|\lambda|=\lambda_1 +\cdots +\lambda_m$
and whose length $\ell(\lambda)=m$.  Each partition $\lambda$ has an 
associated Ferrers diagram with $\lambda_i$ lattice squares in the $i^{th}$ 
row, from the bottom to top (French notation).  Any lattice square $(i,j)$ 
in the $i$th row and $j$th column of a Ferrers diagram is called a cell.
The conjugate of $\lambda$, denoted $\lambda'$, is the reflection of
$\lambda$ about the main diagonal.  
$\lambda$ is ``{\it $k$-bounded}'' if $\lambda_1 \leq k$ and the set of all 
such partitions is denoted $\mathcal P^k$.  The partition
$\lambda\cup\mu$ is the non-decreasing rearrangement of the parts 
of $\lambda$ and $\mu$.  We say that 
$\lambda \subseteq \mu$ when $\lambda_i \leq \mu_i$ for
all $i$.  Dominance order $\unrhd$ is defined on partitions by
$\lambda\unrhd\mu$ when $\lambda_1+\cdots+\lambda_i\geq
\mu_1+\cdots+\mu_i$ for all $i$, and $|\lambda|=|\mu|$.

\medskip

More generally, for $\rho \subseteq \gamma$, the skew shape $\gg/\rho$ is 
identified with its diagram $\{(i,j) : \rho_i<j\leq \gg_i\}$.
Lattice squares that do not lie in $\gg/\rho$ will be called 
``{\it squares}'' instead of cells. 
We say that any $c\in \rho$ lies ``{\it below}'' $\gg/\rho$.
The hook of any lattice square $s\in \gg$ is defined as the 
collection of cells  of $\gg/\rho$ that lie inside the $L$ with $s$ as 
its corner.  This applies to all $s\in \gg$ including those 
below $\gg/\rho$.  For example, the hook of $s=(1,3)$ is  
depicted by the framed cells:
\begin{equation}
\gg/\rho\,=\,(5,5,4,1)/(4,2)\,=\,
{\tiny{\tableau*[scY]
{\cr&&\tf&\cr\bl&\bl&\tf&&\cr\bl &\bl &\bl s &\bl & \tf}}} \, .
\end{equation}
The hook-length of $s$, $h_s(\gg/\rho)$, is the number of 
cells in the hook of $s$.  In the  example, 
$h_{(1,3)}\big((5,5,4,1)/(4,2)\big)=3$ 
and $h_{(3,2)}\big((5,5,4,1)/(4,2)\big)=3$.  
A cell or square has a \dit{$k$-bounded hook} if its hook-length is
no larger than $k$.  For a partition $\lambda$, $h(\lambda)$ 
refers to hook-length of cell (1,1) called the main hook-length.

\medskip

A {$p$-core is a partition that does not contain any hooks of 
length $p$, and $\mathcal C^p$ will denote the set of all $p$-cores.
The $p$-residue of square  $(i,j)$ is $j-i \mod p$;
the label of this square when squares are periodically labeled 
with $0,1,\ldots,p-1$, where zeros lie on the main diagonal
(see \cite{[JK]} for more on cores and residues).  The 5-residues 
associated to the 5-core $(6,4,3,1,1,1)$ are
$$
{\tiny{\tableau[scY]{\bl 4 |0|1|2,\bl 3|3,4,0,\bl 1|4,0,1,2,
\bl 3| 0,1,2,3,4,0,\bl 1}}}
$$
A ``{\it removable}'' corner of
partition $\gg$ is a cell $(i,j)\in \gg$ with $(i ,j+1),(i+1,j)\not\in \gg$
and an ``{\it  addable}'' corner of $\gg$ is a square  
$(i,j)\not\in \gg$ with $(i ,j-1),(i-1,j)\in \gg$.
Note, squares $(\ell(\gg)+1,1)$ and $(1,\gg_1+1)$
are also considered addable.
\begin{remark} \label{noboth}
A given $p$-core never has  a removable and an addable corner of the 
same residue  (e.g.  \cite{[LMcore]}).
\end{remark}

\medskip

A (semi-standard or column-strict) tableau $T$ is a filling of a 
Ferrers shape with integers that strictly increase in columns and weakly 
increase in rows.  The weight of $T$ is the composition $\alpha$ 
where $\alpha_i$ is the multiplicity of $i$ in the tableau $T$.

\section{Definition of $k$-Schur functions} 
\label{sect3}

Proving the beautiful properties that were conjectured to be held by 
the $k$-Schur functions has not come easily with the prior characterizations 
for these functions.  However, as discussed in the introduction,
a lengthy empirical study of this basis led to a family of 
tableaux defined by certain fillings of $k+1$-cores.  These tableaux
connect directly to the type-A affine Weyl group and conjecturally
enumerate the monomial terms in the $k$-Schur expansion of Macdonald 
polynomials \cite{[LMcore]}.  This family of tableaux 
is the central object of study here -- producing a combinatorial definition 
for $k$-Schur functions that enables us to prove properties still
conjectural for the earlier characterizations.

\begin{definition}
\label{defktabgen} 
Let $\gg$ be a $k+1$-core with $m$ $k$-bounded hooks and  let
$\aa=(\aa_1,\ldots,\aa_r)$ be a composition of $m$. 
A ``$k$-tableau" of shape $\gg$ and ``$k$-weight" $\aa$ 
is a filling of $\gg$ with integers $1,2,\ldots,r$ such that

\smallskip
\noindent
(i) rows are weakly increasing and columns are strictly increasing

\smallskip
\noindent
(ii) the collection of cells filled with letter $i$ are labeled by exactly 
$\alpha_i$ distinct $k+1$-residues.
\end{definition}

\begin{example}
\label{exssktab}
The $3$-tableaux of $3$-weight $(1,3,1,2,1,1)$ and shape $(8,5,2,1)$ are:
\begin{equation}
{\tiny{\tableau*[scY]{5\cr 4&6\cr2&3&4&4&6\cr 1&2&2&2&3&4&4&6 }}} \quad
{\tiny{\tableau*[scY]{6\cr 4&5\cr 2&3&4&4&5\cr 1&2&2&2&3&4&4&5 }}} \quad
{\tiny{\tableau*[scY]{4\cr 3&6\cr 2&4&4&5&6\cr 1&2&2&2&4&4&5&6 }}}
\end{equation}
\end{example}

\medskip

\noindent 
More generally a notion of skew $k$-tableaux follows naturally
from $k$-tableaux, see Definition~\ref{defktabgenskew}.

\medskip

\begin{remark}
\label{ktabtab}
When $k\geq h(\gamma)$, a $k$-tableau $T$ of shape $\gamma$ and $k$-weight 
$\mu$ is a semi-standard tableau of weight $\mu$ 
since no two diagonals of $T$ can have the same residue. 
\end{remark}

Although a $k$-tableau is associated to a shape $\gamma$ and weight
$\alpha$, in contrast to usual tableaux, $|\alpha|$ does not equal $|\gamma|$ 
in general. Instead, $|\alpha|$ is the number of $k$-bounded hooks 
in $\gamma$.  This distinction is natural in view of a bijective
correspondence between $k+1$-cores and  $k$-bounded
partitions that was defined in \cite{[LMcore]} by the map:
$$
\kbnd: \CC^{k+1}\to \CP^k\quad
\text{where}\quad 
\kbnd\left(\gg\right) = (\lambda_1,\ldots,\lambda_\ell)
\,,
$$
with $\lambda_i$ denoting the number of cells with a $k$-bounded hook
in row $i$ of $\gamma$.  Note that the number of $k$-bounded hooks 
in $\gamma$ is $|\lambda|$.  The inverse map relies on constructing a 
certain skew diagram $\gg/\rho$ from 
$\lambda$, and setting $\core(\lambda)=\gamma$.
These special skew diagrams are defined:

\begin{definition} \label{def5}
For $\lambda\in\mathcal P^k$, the
``$k$-skew diagram of $\lambda$" is the diagram $\lambda/^k$
where

(i) row $i$ has length $\lambda_i$ for $i=1,\ldots,\ell(\lambda)$

(ii) no cell of $\lambda/^k$ has hook-length exceeding $k$

(iii) all squares below $\lambda/^k$ have hook-length exceeding $k$.

\noindent
\end{definition}

A convenient algorithm for constructing the diagram of $\lambda/^k$ is 
given by successively attaching a row of length $\lambda_i$ to the bottom of 
$(\lambda_1,\dots,\lambda_{i-1})/^k$ in the leftmost position so that no 
hook-lengths exceeding $k$ are created.

\begin{example} \label{exskew}
Given $\lambda =(4,3,2,2,1,1)$ and $k=4$,
\begin{equation*}
\lambda = {\tiny{\tableau*[scY]{  \cr  \cr  & \cr &
\cr & & \cr & & & }}}
\quad
\implies
\quad
\lambda/^4 = {\tiny{\tableau*[scY]{  
\cr  
\cr  & 
\cr \bl &  & 
\cr  \bl &\bl & & & 
\cr \bl &\bl& \bl & \bl &\bl & & & &
}}}
\qquad\implies
\;
\core(\lambda)= {\tiny{\tableau*[scY]{  
\cr  
\cr  & 
\cr  &  & 
\cr  & & & & 
\cr  && & & & & & &
}}}
\end{equation*}
\end{example}

\medskip

The bijection between $k+1$-cores and $k$-bounded partitions gives
rise to a natural {\it involution} on the set of $k$-bounded partitions
that refines partition conjugation.

\begin{definition}
\cite{[LMcore]}
The \dit{$k$-conjugate} of a $k$-bounded partition $\lambda$ is
$$
\lambda^{\omega_k} := \core^{-1}\left(\core(\lambda)'\right)
$$
\end{definition}

\begin{remark} 
\label{conjupetit}
$\lambda^{\omega_k}=\lambda'$  when $h(\lambda) \leq k$
since $\core(\lambda)=\lambda$ in this case.
\end{remark}

\medskip

The analogy with usual tableaux is now more apparent.  We denote the 
set of all $k$-tableaux of shape $\core(\mu)$ and $k$-weight $\aa$
by $\mathcal T^k_{\aa}(\mu)$, and call a $k$-tableau \dit{standard} 
when its $k$-weight is $(1^n)$.  Here, we will study 
properties of the ``$k$-Kostka numbers":
\begin{equation}
K_{\mu\alpha}^{(k)}\;=\;|T^k_\aa(\mu)|
\,.
\end{equation}
For example, they
satisfy a triangularity property similar to the Kostka numbers. 

\begin{property}
\label{trikostka}
\cite{[LMcore]} 
For any $k$-bounded partitions $\lambda$ and $\mu$,
\begin{equation}
K_{\mu\lambda}^{(k)}=0 \quad\text{when}\quad \mu \ntrianglerighteq
\lambda \quad
\text{ and }\quad K_{\mu\mu}^{(k)}=1\,.
\end{equation}
\end{property}

Thus the matrix $||K^{(k)}||_{\lambda,\mu\in \mathcal P^k}$ is invertible, 
naturally giving rise to a family of functions defined by:

\begin{definition}
\cite{[LMcore]}
\label{kschurdef}
The ``$k$-Schur functions", indexed by $k$-bounded partitions,
are defined by inverting the unitriangular system:
\begin{equation}
\label{e1}
h_\lambda = s_\lambda^{(k)}+\sum_{\mu : \mu\rhd\lambda} K_{\mu\lambda}^{(k)}
s_\mu^{(k)}\,\quad\text{for all } \lambda_1\leq k\,.
\end{equation}
\end{definition}

\medskip

This characterization has the advantage that properties of $k$-Schur 
functions can be derived from the study of $k$-tableaux.  
The next sections will thus be devoted to understanding the $k$-tableaux
and the proofs of $k$-Schur function properties will follow in 
\S\ref{pierisec}.

\section{Properties of $k$-tableaux}
\label{ktabprops}

Since the homogeneous symmetric functions commute, Definition~\ref{kschurdef} 
suggests that $K^{(k)}_{\lambda\alpha}=K^{(k)}_{\lambda\mu}$ for any 
rearrangement $\mu$ of the parts of $\alpha$.  In \cite{[LMcore]},
it was conjectured that this could be explained by finding an involution 
on the set of $k$-tableaux sending $\mathcal T^k_\alpha(\lambda)$ to 
$\mathcal T^k_{\hat\alpha}(\lambda)$, where $\hat\alpha$ is obtained 
by transposing two adjacent components of $\alpha$.  The construction 
of such an involution is the spring board to proving properties of the 
$k$-Schur functions. To this end, we now explore characteristics of 
$k$-tableaux that are necessary to construct and prove the involution.

\medskip

We first earmark certain cells of a partition $\lambda$.
A cell $(i,j)$ where $(i+1,j+1)\not\in \la$ is called \dit{extremal}.
The squares $(0,\lambda_1)$ and $(\ell(\lambda),0)$, those below 
and to the left of the diagram of $\lambda$, will also be called 
extremal.  The cell immediately to the left (or right) of the cell $c$ 
is {\it ``left-adj"} (or {\it ``right-adj"}) to $c$, while the 
cell $(i-1,j-1)$ is {\it ``south-west"} of $(i,j)$. 
We will repeatedly use the following property of 
cores:

\begin{remark}
\label{core}
In a $p$-core $\gamma$, if an extremal lies at the top of its column 
in some row $r$, then in all rows weakly lower than $r$, extremal of 
the same residue must lie at the top of their column.  Similarly, if an 
extremal lies at the end of its row in some row $r$, then in 
all rows weakly higher than $r$, all extremals of the same residue 
lie at the end of their row.  Note that this argument applies to all 
extremals, including those that are not cells -- $(0,\gamma_1)$ or 
$(\ell(\gamma),0)$.
\end{remark}

With the goal of producing an involution that switches the
weight of consecutive letters $a$ and $b=a+1$ in a $k$-tableau, the 
behavior of these letters is our main concern.  We consider
entries $a$ and $b$ to be {\it ``married"} if they occur in the 
same column, an entry $a$ (resp. $b$) is {\it ``divorced"} if it has 
the same residue as some married $a$ (resp. $b$), and {\it ``single"} 
otherwise.  When the letter $x$ occupies a cell in row $r$ that is labeled with 
residue $j$, we say this cell contains an $x_r(j)$, or simply an $x(j)$.  
$Res_r(x)$ will be the set of all residues that label cells occupied by a 
letter $x$ in row $r$, while $URes_r(x)$ will be only the residues labeling 
unmarried $x$'s in row $r$.  We also consider 
$URes_{r}(a,b)=URes_{r}(a)\cup URes_{r}(b)$.  

\medskip

The Bender-Knuth involution for semi-standard tableau is based on 
the following simple observation, also needed for our purposes:

\begin{remark}
\label{rem}
The married $b$'s in row $r$ lie at the end of the sequence of $b$'s in 
that row since an unmarried $b$ must have an entry smaller than
$a$ lying below it.  Similarly, married $a$'s in any given row lie 
at the beginning of the sequence of $a$'s in that row.
\end{remark}

The extension of their involution to $k$-tableaux requires several
intricate properties whose proofs rely heavily on the fact 
that deleting a letter from a $k$-tableau gives a $k$-tableau.  
To be precise, Proposition~32 in \cite{[LMcore]} states that deleting 
the largest letter from a standard $k$-tableau produces a new 
standard $k$-tableau, and thus by iteration, deleting all letters 
larger than any letter $x$ produces a $k$-tableau.   This idea
can be generalized to the non-standard $k$-tableau case 
by introducing a total order on the letters/residues. 
In particular, consider the list of residues $(j_1,\ldots, j_r)$ occupied 
by $x$, where $j_1$ is the residue of the lowest, rightmost $x$ in $T$, 
$j_2$ is the residue of the lowest, rightmost $x$ in $\hat T$ (obtained 
by deleting $x(j_1)$ from $T$), and so forth.  It was shown in \S 10.1 
of \cite{[LMcore]} that replacing $x(j_1)$ with $n$, $x(j_2)$ with $n-1$, 
and so forth produces a standard $k$-tableau.  Therefore, using the 
order $z>x(j_1)>x(j_2)>\cdots>x(j_r)>y$ for $z>x>y$, a 
$k$-tableau $T_{\leq x(i)}$ 
(resp. $T_{< x(i)}$) is obtained by deleting all letters larger 
(weakly larger) than $x(i)$ from $T$.

\begin{remark}
\label{remorder}
It is important to note that $x(i)<x(i+1)$ when both $x(i),x(i+1)$ are 
in a $k$-tableau $T$.  Otherwise, $x(i)>x(i+1)$ would imply that 
$x(i),x(i+1)\in T_{\leq x(i)}$, where $T_{\leq x(i)}$ is a $k$-tableau 
with an extremal $x_r(i)$ at the end of some row $r$ and 
an extremal of residue $i$ left-adj to an $x_m(i+1)$ for $m>r$,
contradicting Remark~\ref{core}.
\end{remark}

\medskip

\begin{property}
\label{tfae} $~$
\begin{enumerate}
\item[(i)]
Given an unmarried $b(j)$ in a $k$-tableau $T$, any $a(j)\in T$ 
is married and lies weakly higher than the highest unmarried $b(j)$.
Further, $a(j)$ occurs in $T$ if and only if there is a divorced 
$b(j-1)$ left-adj to the unmarried $b(j)$.  

\smallskip

\item[(ii)]
Given an unmarried $a(j)$ in a $k$-tableau $T$, any $b(j)\in T$ 
is married and lies strictly higher than the highest unmarried $a(j)$.
Further, $b(j)$ occurs in $T$ if and only if there is a divorced 
$a(j+1)$ right-adj to the unmarried $a(j)$.  
\end{enumerate}
\end{property}
\begin{proof}
We prove case (i) and note that the other case follows similarly.
Given $a_r(j)\in T$, we first claim it does not lie lower than any
unmarried $b(j)$.  Suppose there is an unmarried $b_m(j)$ for some 
$m>r$.  In $T_{\leq a(j)}$, $a_r(j)$ lies at the end of its row by 
Remark~\ref{remorder} implying that all extremals above row $r$ lie 
at the end of their row by Remark~\ref{core}.  However,
the extremal south-west of $b_m(j)$ lies at the end of its 
row in $T_{\leq a(j)}$ only if an $a(j+1)$ lies below the 
$b_m(j)$ in $T$, contradicting that $b_m(j)$ is unmarried.
Thus, given an unmarried $b_m(j)\in T$ and an $a_r(j)$ with
$r\geq m$, it remains to show that $a_r(j)$ is married.  In 
$T_{<b(j)}$, the cell of residue $j-1$ left-adj to $b_m(j)$ 
lies at the end of its row.  Thus since $a_r(j)\in T_{<b(j)}$,
$a_r(j)$ must be married to prevent its left-adj cell of residue 
$j-1$ from being extremal in $T_{<b(j)}$.

For the second part of the assertion, the $\Leftarrow$ implication 
holds since a divorced $b(j-1)\in T$ implies there is a married 
$b(j-1)\in T$, lying above an $a(j)$.  For the $\Rightarrow$ implication, 
consider $T$ with an unmarried $b_r(j)$ and some $a_m(j)$.  The previous 
paragraph explains that $m\geq r$ and $a_m(j)$ is married.  
In $T_{<b(j-1)}$, since $a_m(j)$ lies at the top of its column,
the extremal south-west of $b_r(j)$ lies at the top of its column.  
The only way the cell left-adj to $b_r(j)$ is not in $T_{<b(j-1)}$ is
if it is $b_r(j-1)$.  This $b_r(j-1)$ is not single since $a_m(j)$ is 
married to $b(j-1)$. Further, it is not married since no $a(j)$ lies 
lower than row $m\geq r$.
\end{proof}

\begin{lemma}
\label{eq7}
In a $k$-tableau $T$ with an $x_r(i)$ and an $x_m(i)$ for $r<m$, 
$Res_m(x)\subseteq Res_r(x)$.
\end{lemma}
\begin{proof}
If an $x_m(i+1)\in T$, then it lies at the top of its column in 
$T_{\leq x(i+1)}$.  Thus, there must be an $x_r(i+1)$ right-adj 
to $x_r(i)$ to prevent the entry of residue $i+1$ below $x_r(i)$ 
from being extremal in $T_{\leq x(i+1)}$.  If an $x_{m}(i-1)\in T$, 
then the cell of residue $i$ below $x_m(i-1)$ lies at the top of its
column in $T_{<x(i-1)}$.  Therefore, an $x_r(i-1)$ must be left-adj to 
$x_r(i)$ to ensure that the extremal south-west of $x_r(i)$ lies 
at the top of its column in $T_{<x(i-1)}$.  By iteration, 
$Res_m(x)\subseteq Res_r(x)$.  
\end{proof}

\begin{property}\label{lemmaint}
Let $URes_{m}(a,b)\cap URes_{r}(a,b) \neq \emptyset$ for some $r<m$.  Then
$$
URes_m(a)\subseteq URes_r(a)\quad\text{and}\quad
URes_m(b)\subseteq URes_r(b)\,.
$$
In particular, this implies $URes_m(a,b)\subseteq URes_r(a,b)$. 
\end{property}
\begin{proof} 
Since an unmarried $a$ and $b$ of the same residue cannot lie in $T$ 
by Property~\ref{tfae}, $URes_{m}(a,b)\cap URes_{r}(a,b)\neq\emptyset$ 
implies $URes_{m}(a)\cap URes_{r}(a)\neq\emptyset$ or 
$URes_{m}(b)\cap URes_{r}(b)\neq\emptyset$.  Thus
$Res_m(a)\subseteq Res_r(a)$ or $Res_m(b)\subseteq Res_r(b)$
by Lemma~\ref{eq7}.  Note that $Res_m(a)\subseteq Res_r(a)$ 
implies $URes_m(a)\subseteq URes_r(a)$ since
an unmarried $a_m(i)\in T$ lies at the top of its column in
$T_{\leq b(i-1)}$ forcing $a_r(i)$ also to lie 
at the top of its column. Similarly we have
$URes_m(b)\subseteq URes_r(b)$.  Therefore, $URes_m(a)\subseteq URes_r(a)$ 
{\bf or} $URes_m(b)\subseteq URes_r(b)$ and it remains to show that 
both in fact are true.  

It suffices to prove that $\emptyset \neq URes_m(a)\subseteq URes_r(a)$ 
implies $URes_m(b)\subseteq URes_r(b)$ and to note by a similar argument
that $\emptyset \neq URes_m(b)\subseteq URes_r(b)$
implies $URes_m(a)\subseteq URes_r(a)$.
Let $a_m(i)$ denote the rightmost $a$ in row $m$.  If any unmarried
$b$ lies in row $m$, then there is an unmarried $b_m(i+1)$ right-adj 
to $a_m(i)$ by Remark~\ref{rem}.  Since $b_m(i+1)$ lies at the top of 
its column in $T_{\leq b(i+1)}$, there must be an entry $x\leq b(i+1)$ 
right-adj to $a_r(i)$ to prevent the entry of residue $i+1$ below 
$a_r(i)$ from being extremal in $T_{\leq b(i+1)}$.  Property~\ref{tfae}(i) 
ensures that $x\neq a_r(i+1)$ since there is an unmarried $b_m(i+1)$.  
Therefore $x=b_r(i+1)\in T$.  Thus we have $b_m(i+1),b_r(i+1)\in T$
and can use Lemma~\ref{eq7} to obtain $Res_m(b)\subseteq Res_r(b)$.
\end{proof}

Our involution will be defined on certain rows of a $k$-tableau
that are characterized by the following equivalence relation.

\begin{definition} 
Rows $r_1$ and $r_2$ in a $k$-tableau are equivalent, 
$``\,r_1 \sim_a r_2\,"$, when they satisfy the following conditions:
\begin{itemize}
\item $URes_{r_1}(a,b)\neq \emptyset$ and $URes_{r_2}(a,b)\neq \emptyset$ 
\item $URes_{r_1}(a,b) \subseteq URes_{r}(a,b)$ and 
$URes_{r_2}(a,b) \subseteq URes_{r}(a,b)$ for some $r$.
\end{itemize}
\end{definition}

\begin{proposition} 
$\sim_a$ is an equivalence relation on the set of
rows in a $k$-tableau containing an unmarried $a$ or $b$.
\end{proposition}
\begin{proof}  
The only non-trivial part is to show transitivity.  
With $r_1 \sim_a r_2$ and $r_2 \sim_a r_3$, 
\begin{eqnarray*}
\emptyset \neq  URes_{r_1}(a,b) \subseteq URes_{r}(a,b)\;\;\text{and}\;\;
\emptyset \neq  URes_{r_2}(a,b) \subseteq URes_{r}(a,b)\;\;
\text{for some}\;\; r\\
\emptyset \neq  URes_{r_2}(a,b) \subseteq URes_{t}(a,b)\;\;\text{and}\;\;
\emptyset \neq  URes_{r_3}(a,b) \subseteq URes_{t}(a,b)\;\;
\text{for some}\;\; t\,.
\end{eqnarray*}
Thus, in particular, $URes_{r}(a,b) \cap URes_{t}(a,b)\neq \emptyset$
giving that $URes_{\max(r,t)}(a,b) \subseteq URes_{\min(r,t)}(a,b)$
by Property~\ref{lemmaint}.  Therefore 
$$URes_{r_1}(a,b) \subseteq URes_{\min(r,t)}(a,b) \quad
{\rm and} \quad URes_{r_3}(a,b) \subseteq URes_{\min(r,t)}(a,b)\, ,$$
implying that $r_1 \sim_a r_3$.
\end{proof}

We can take the lowest row in each equivalence class for a set of
representatives.  Property~\ref{lemmaint} implies that these
representatives can equivalently be defined by:

\begin{definition} A ``\,fundamental row" of a tableau is a row $m$ where 
$URes_m(a,b)$ is not contained in $URes_r(a,b)$ for any $r<m$.
\end{definition}

\section{An involution on $k$-tableaux}

We are now ready to construct an involution on the set of $k$-tableaux 
sending $\mathcal T^k_\alpha(\lambda)$ to $\mathcal T^k_{\hat\alpha}(\lambda)$, 
where $\hat\alpha$ is obtained by transposing two adjacent components 
of $\alpha$.   Recall that the Bender-Knuth involution \cite{[BK]} is 
defined on semi-standard tableau by sending the string $a^rb^s$ of 
single $a$'s and $b$'s in each row to the
string $a^s b^r$, thus permuting the weight of the tableau.  In our 
case, we perform a similar operation but must take into consideration 
the added notion of divorced entries.  Our algorithm boils down to
applying the BK involution to fundamental rows and ``correcting".
It reduces to the BK involution for large $k$.

\begin{definition}
\label{definv}
The operator $\tau_a$ on a $k$-tableau $T$ is defined as follows on
the equivalence classes $\mathcal C_i=\{r \, |\, r \sim_a r_i \}$, 
for the set of fundamental rows $r_1,\ldots,r_n\in T$:
\begin{enumerate}
\item In row $r_i$:
\begin{enumerate}
\item Replace the entries $a^tb^s$ of \underline{\it single} 
$a$'s and $b$'s by $a^sb^t$.

\smallskip
\item If $t>s$, relabel any $a$ lying to the right of some $b$ by a $b$.
Otherwise, relabel any $b$ lying to the left of an $a$ with an $a$.

\end{enumerate}

\smallskip

\item  In rows above $r_i$: for $\mathcal S_i$ the set of 
residues of $a$'s (or $b$'s) that were relabeled in step 1,
correspondingly relabel every unmarried $a$ (or $b$) 
that has residue in $\mathcal S_i$. 
\end{enumerate}
\end{definition}

Note by definition of $\sim_a$ that step 2 only involves 
rows in the class $\mathcal C_i$ implying that no row is involved
in this step for two distinct values of $i$.

\begin{example}
Given a $4$-tableau of weight (2,1,4,2,3), we act with $\tau_4$ 
to permute the number of residues occupied by letters 4 and 5.
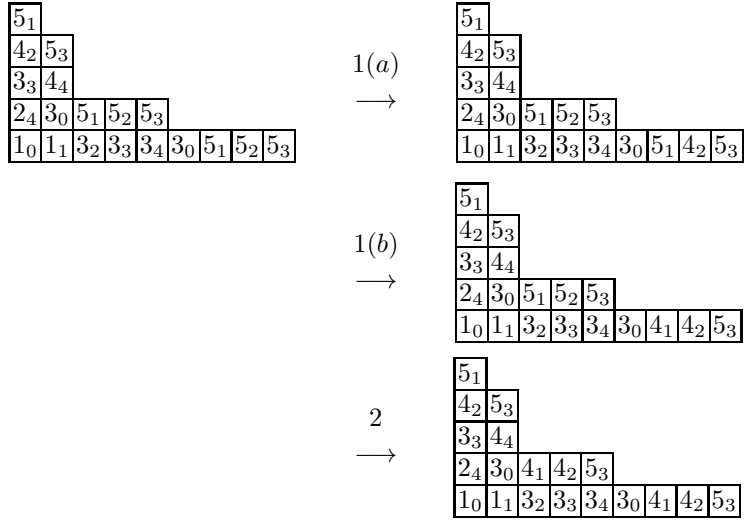
\begin{equation*}
{{\tableau*[scY]{
5_1 \cr
4_2 & 5_3 \cr
3_3 & 4_4 \cr
2_4 & 3_0 & 5_1 & 5_2 & 5_3 \cr
1_0 & 1_1 &  3_2 & 3_3 & 3_4 & 3_0 & 5_1 & 5_2 & 5_3 }}}
\qquad
{1(a) \atop
\longrightarrow }
\qquad
{{\tableau*[scY]{
5_1 \cr
4_2 & 5_3 \cr
3_3 & 4_4 \cr
2_4 & 3_0 & 5_1 & 5_2 & 5_3 \cr
1_0 & 1_1 &  3_2 & 3_3 & 3_4 & 3_0 & 5_1 & 4_2 & 5_3 }}}
\end{equation*}
\begin{equation*}
\phantom{{\tableau*[scY]{
5_1 \cr
4_2 & 5_3 \cr
3_3 & 4_4 \cr
2_4 & 3_0 & 5_1 & 5_2 & 5_3 \cr
1_0 & 1_1 &  3_2 & 3_3 & 3_4 & 3_0 & 5_1 & 4_2 & 5_3 }}}
\qquad
{1(b) \atop
\longrightarrow }
\qquad
{{\tableau*[scY]{
5_1 \cr
4_2 & 5_3 \cr
3_3 & 4_4 \cr
2_4 & 3_0 & 5_1 & 5_2 & 5_3 \cr
1_0 & 1_1 &  3_2 & 3_3 & 3_4 & 3_0 & 4_1 & 4_2 & 5_3 }}}
\end{equation*}
\begin{equation*}
\phantom{{\tableau*[scY]{
5_1 \cr
4_2 & 5_3 \cr
3_3 & 4_4 \cr
2_4 & 3_0 & 5_1 & 5_2 & 5_3 \cr
1_0 & 1_1 &  3_2 & 3_3 & 3_4 & 3_0 & 5_1 & 4_2 & 5_3 }}}
\qquad
{2 \atop
\longrightarrow }
\qquad
{{\tableau*[scY]{
5_1 \cr
4_2 & 5_3 \cr
3_3 & 4_4 \cr
2_4 & 3_0 & 4_1 & 4_2 & 5_3 \cr
1_0 & 1_1 &  3_2 & 3_3 & 3_4 & 3_0 & 4_1 & 4_2 & 5_3 }}}
\end{equation*}
\end{example}

\medskip

Although it is not immediately clear that the number of
residues occupied by 4's and 5's in the $k$-tableaux
has been switched, we will use properties from the previous 
section to prove that $\tau_a$ does in fact change
$k$-tableaux as desired.

\medskip

\begin{proposition}
For any $T \in \mathcal T^k_\alpha(\lambda)$, the tableau 
$\tau_a (T)$ belongs to $\mathcal T^k_{\hat\alpha}(\lambda)$,  
where $\hat \alpha=(\ldots,\alpha_{a+1},\alpha_a,\ldots,)$
is obtained by transposing $\alpha_{a}$ and $\alpha_{a+1}$ in $\alpha$.
\end{proposition}
\begin{proof}
We start by showing that $\hat T=\tau_a (T)$ is a column-strict
tableau. Then proving that the weight changes 
in the specified manner will imply it is a $k$-tableau.
By Remark~\ref{rem}, $T$ has a non-decreasing contiguous sequence of 
unmarried letters $a$ and $b$.  Thus in Step 1, only unmarried
$a$ and $b$ in rows $r_i$ are changed, and the definition of $\tau_a$ 
implies that these rows of $\hat T$ are non-decreasing.  In Step 2, 
rows of $T$ in the classes $\mathcal C_i$ are changed according to 
entries that were changed in Step 1.   Since the unmarried $a$'s and $b$'s 
in such rows form a (contiguous) subsequence of the unmarried letters 
in row $r_i$ by Property~\ref{lemmaint}, these rows are also non-decreasing 
in $\hat T$.  Further, since an unmarried $a$ lies below an entry strictly 
larger than $b$ while an unmarried $b$ lies above an entry strictly smaller 
than $a$, changing an unmarried $a$ to $b$ or $b$ to $a$ retains the
property of strictly increasing columns.

This given, it remains to show that $\hat\alpha_a=\alpha_b$ and 
$\hat\alpha_b=\alpha_a$.  Let $\alpha_a^s$ and $\alpha_a^m$ denote the
number of single (resp. married) residues occupied by letter $a$ in $T$,
and observe that $\alpha_a = \alpha_a^s+\alpha_a^m$.   
Similarly for  the letter $b$.  Since a married $a$ or $b$ remains as 
such under the action of $\tau_a$, we have that
$\hat\alpha_a^m=\alpha_a^m=\alpha_b^m=\hat\alpha_b^m$.
Thus, we need only show that
$\hat\alpha_a^s = \alpha_b^s$ and $\hat\alpha_b^s=\alpha_a^s$.
However, by Property~\ref{lemmaint},
\begin{equation}
\alpha_a^s=\sum_i\alpha_a^s(r_i) 
\quad\text{and}\quad
\alpha_b^s=\sum\alpha_b^s(r_i)
\,,
\end{equation}
for $\alpha_a^s(r_i)$ the number of single $a$ residues in row $r_i$ 
of $T$.   Since the definition of $\tau_a$ implies that
each single $a(j)$ or $b(j)$ occurs in exactly one 
fundamental row $r_i$ of $\hat T$, we can further reduce our problem to showing
\begin{equation}
\hat\alpha_a^s(r_i)=\alpha_b^s(r_i)
\quad\text{ and }\quad\hat\alpha_b^s(r_i)=\alpha_a^s(r_i)\,.
\end{equation}

To show that the number of single $a$ and $b$ residues is permuted in 
a fundamental row $r_i$,  first note in Step 1(a) that when single 
$b$'s are relabeled by $a$'s, the number of single $b$-residues is 
decreased by $\alpha^s_b(r_i)-\alpha^s_a(r_i)$ (considering the 
case $\alpha^s_b(r_i)>\alpha^s_a(r_i)$).  Since Step 1(b) involves 
only divorced entries, no further single $b$-residues are lost
implying that $\hat\alpha_b^s(r_i)= \alpha_a^s(r_i)$.  

To prove $\hat\alpha_a^s(r_i)= \alpha_b^s(r_i)$, we must
verify that precisely $\alpha^s_b(r_i)-\alpha^s_a(r_i)$ $b$'s are 
sent to {\it single} $a$'s.  Each relabeled $b$ goes to either a 
single or divorced $a\in\hat T$.  To be precise, a $b(j)$ 
goes to a divorced $a(j)$ only if there is an $a(j)\in T$,
and Property~\ref{tfae}(i) tells us $a(j)\in T$ iff
there is a divorced $b(j-1)$ left-adj to $b(j)$. 
Therefore, of the $\alpha^s_b(r_i)-\alpha^s_a(r_i)$ 
$b$'s relabeled in Step 1(a), each $b(j)$ that is right-adj 
to a divorced $b(j-1)$ does not give rise to a single $a$.
However, each of these $b(j)$'s can be matched with
a $b(i)$ that is not right-adj to a divorced $b(i-1)$ and
thus goes to a single $a$ in Step 1(b).
For example,
\begin{eqnarray}
\begin{matrix}
 b_s/a& b_d & b_d & b_d & b_d & b_s\quad & b_s & b_s/b\\
 \downarrow & & &  &\quad \swarrow& & \downarrow&\downarrow\\
 \downarrow & & & \swarrow & & & \downarrow&\downarrow\\
 \downarrow & &\swarrow\quad & & & & \downarrow&\downarrow\\
a_s/a& \quad a_s & a_d & a_d & a_d & a_d & a_s & a_s/b
\end{matrix}
\end{eqnarray}
Therefore, exactly $\alpha^s_b(r_i)-\alpha^s_a(r_i)$ new
single $a$ residues arise and we have proven our assertion.
The case $\alpha^s_a(r_i)>\alpha^s_b(r_i)$ is similar.
\end{proof}

\medskip

\begin{proposition} 
The operator 
$\tau_a$ is an involution on $\mathcal T^k_\alpha(\lambda)$, 
for all $1\leq a<\ell(\alpha)$.
\end{proposition}
\begin{proof}
Since $\tau_a$ acts by sending certain unmarried $a$ to 
unmarried $b$ or vice-versa, the sets $URes_i(a,b)$ are fixed
by $\tau_a$.  In particular, the fundamental rows of 
$T$ are those of $\hat T$.  Property~\ref{tfae}(ii) and an 
argument similar to the one given in the previous proposition
implies further that
$\bigl(\tau_a\bigr)^2$ fixes the entries in the 
fundamental rows of $T$.
For example, applying $\tau_a$ to our 
previous illustration:
\begin{eqnarray}
\begin{matrix}
a_s/a& a_s\quad & a_d & a_d & a_d & a_d & a_s & a_s/b\\
 \downarrow &\quad \searrow & &  && & \downarrow&\downarrow\\
 \downarrow & & & \searrow  & & & \downarrow&\downarrow\\
 \downarrow & & &&\quad\searrow & & \downarrow&\downarrow\\
 b_s/a& b_d & b_d & b_d & b_d & b_s& b_s& b_s/b\\
\end{matrix}
\end{eqnarray}
By definition of $\tau_a$, the equivalence classes 
$\mathcal C_i = \{r \, | \, r\sim_a r_i \}$ are then also
left unchanged by $\bigl(\tau_a\bigr)^2$ and the claim follows.
\end{proof}

The two previous propositions immediately imply:

\begin{theorem} 
\label{bigkostka}
Given $\lambda\in\mathcal P^k$, $\alpha$ a composition of $|\lambda|$,
and any $1\leq a<\ell(\alpha)$, 
$$
\tau_a : \mathcal T^k_\alpha(\lambda)\to
\mathcal T^k_{\hat\alpha}(\lambda)$$
is a bijection, where
$\hat \alpha=(\ldots,\alpha_{a+1},\alpha_a,\ldots,)$.
\end{theorem}
Given that $K^{(k)}_{\lambda\alpha} = | \mathcal T^k_\alpha(\lambda)|$,
the theorem has the following corollary.
\begin{corollary} 
\label{corosym}
For $\lambda\in \mathcal P^k$, and a composition $\alpha$ of $|\lambda|$,
\begin{equation}
K^{(k)}_{\lambda\alpha}=K^{(k)}_{\lambda\mu}\,,
\end{equation}
where $\mu$ is the weakly decreasing rearrangement of $\alpha$.
\end{corollary}

\medskip

We are also able to derive a recursive formula for the $k$-Kostka 
numbers using a correspondence between $k$-tableaux and certain chains 
of partitions.  Following the notation of \cite{[LMcore]}, 
$\mu,\nu$ are ``{\it $\ell$-admissible}'' when $\mu/\nu$ 
and  $\mu^{\om_k}/\nu^{\om_k}$ are respectively horizontal and vertical 
$\ell$-strips.  More generally,  for any composition
$\aa=(\aa_1,\dots,\aa_r)$, a sequence of partitions
$\left(\mu^{(0)},\mu^{(1)},\cdots,\mu^{(r)}\right)$
is ``$\alpha$-admissible" if $ \mu^{(j)},\mu^{(j-1)}$ are
$\alpha_j$-admissible for all $j$.  A bijection was established 
between the sets: 
\begin{equation}
\label{TDbij}
\mathcal T_\alpha^{k}(\mu)\;\leftrightarrow
\mathcal D^k_{\aa}(\mu):=\left\{ 
(\emptyset=\mu^{(0)},\ldots,\mu^{(r)}=\mu) 
\;\;\text{that are}\;\; \alpha\text{-admissible}
\right\}\,,
\end{equation}
implying in particular that
\begin{equation}
\label{remark23}
\text{
$\mu,\nu$ are $\ell$-admissible iff $\core(\mu)/\core(\nu)$=horizontal 
strip with $\ell$ distinct residues.}
\end{equation}
See \cite{[LMcore]} for the construction and details 
of this correspondence.

\medskip

\begin{corollary}
\label{kostka2}
For $k$-bounded partitions $\mu$ and $\lambda$,
and $0<r\leq k$,
\begin{equation}
\label{kids}
K_{\mu\, {(\ell, \lambda)}}^{(k)} \,=\,
\sum_{{\mu/\nu=\text{horizontal $\ell$-strip}}\atop
\mu^{\omega_k}/\nu^{\omega_k}=\text{vertical $\ell$-strip}}
K_{\nu\lambda}^{(k)}
\,.
\end{equation}
\end{corollary}

\begin{proof}
Since every sequence $(\emptyset=\mu^{(0)},\ldots,\mu^{(r)}=\mu)\in 
\mathcal D^k_{\aa}(\mu)$ has the property that
$\mu,\mu^{(r-1)}$ are $\alpha_r$-admissible,
the cardinality of $\mathcal D^k_{\aa}(\mu)$ satisfies 
the recursion:
$$
| \mathcal D^k_{\aa}(\mu)|= \sum_{ 
\nu: \mu,\nu~are~\alpha_r-admissible}
|\mathcal D^k_{(\aa_1,\dots,\aa_{r-1})}(\nu)|
\,.
$$
The bijection \eqref{TDbij} implies that
$|\mathcal D^k_{\aa}(\mu)|=K_{\mu\alpha}^{(k)}$
for all $\alpha$,  and thus by Corollary~\ref{corosym},
$|\mathcal D^k_{\aa}(\mu)|=K_{\mu\nu}^{(k)}$ for $\nu$ any 
rearrangement of $\alpha$. Therefore,
$$
K_{\mu \,\,  (\ell,\lambda)}^{(k)} 
=|\mathcal D^k_{(\lambda,\ell)}(\mu)|
=
\sum_{ \nu: \mu,\nu~are~\ell-admissible}
|\mathcal D^k_{\lambda}(\nu)| =
\sum_{\nu: \mu,\nu~are~\ell-admissible}
K_{\nu \lambda}^{(k)}
\,.
$$
\end{proof}

\section{Properties of $k$-Schur functions}
\label{pierisec}

In the introduction, we discussed that the $k$-Schur functions have 
been thought to play the role of Schur functions in the spaces spanned 
by homogenous symmetric functions indexed by $k$-bounded partitions:
\begin{equation}
\Lambda^{(k)}=\mathcal L\left\{h_\lambda\right\}_{\lambda_1\leq k}
\,.
\end{equation}
This belief was supported by strong computational evidence that the 
$k$-Schur functions obey refinements of the combinatorial properties 
of Schur functions.  We can now capitalize on our knowledge of 
$k$-tableaux to prove that such beautiful combinatorial properties 
are held by the $k$-Schur functions introduced in
Definition~\ref{kschurdef} by inverting the system
\begin{equation}
h_\lambda = s_\lambda^{(k)}+\sum_{\mu : \mu\rhd\lambda} K_{\mu\lambda}^{(k)}
s_\mu^{(k)}\,\quad\text{for all } \lambda_1\leq k\,.
\end{equation}

\medskip

Immediate from the definition, we have that

\begin{property}
The set $\left\{s_\lambda^{(k)}\right\}_{\lambda_1\leq k}$ forms a 
basis for $\Lambda^{(k)}$.
\end{property}

The unitriangular expression for $h_\lambda$ in terms of $k$-Schurs,
as well as the unitriangular relation between the usual
Schur functions and $h_\lambda$ imply that

\begin{property}
For any $k$-bounded partition $\lambda$,
\begin{equation}
\label{triom}
s_{\lambda}^{(k)}=s_{\lambda}+ \sum_{\mu:\mu\rhd\lambda} d_{\lambda\mu}^{(k)}\, 
s_\mu \;\;\; for\;\; d_{\lambda\mu}^{(k)}\in\mathbb Z\,.
\end{equation}
\end{property}

Although at this point, we can only prove that the coefficients 
$d_{\lambda\mu}^{(k)}$ are integral, we believe that they
are in fact positive.  This would follow by proving that the
$k$-Schur functions studied here are precisely the atoms of 
\cite{[LLM]}, since the atoms are positive sums of Schur functions
by definition.

\medskip

\subsection{Pieri rules}
Much of our prior work with the $k$-Schur functions was drawn from a
conjecture that the $k$-Schur functions satisfy a refinement of the 
Pieri rule called the \dit{$k$-Pieri rule}.  In fact, the characterization 
of the $k$-Schurs used in this article was motivated purely so that they 
would satisfy this rule.  

\medskip

\begin{theorem}
\label{kpieri}
For any $k$-bounded partition $\nu$ and $\ell\leq k$,
\begin{equation}
\label{kpierieq}
h_\ell\,s_\nu^{(k)} = \sum_{\mu\in H_{\nu,\ell}^{(k)}}
s_\mu^{(k)}
\end{equation}
where
$
H_{\nu,\ell}^{(k)}=
\left\{ \mu\, \Big| \, 
\mu/\nu={\rm horizontal~} \ell\text{-}{\rm strip} \;\; {\rm and} \;\;
\mu^{\omega_k}/\nu^{\omega_k}={\rm vertical~} \ell\text{-}{\rm strip} 
\right\} \,.  $
\end{theorem}

\begin{proof}
Since the $k$-Schur functions form a basis of $\Lambda^{(k)}$, 
there is an expansion
\begin{equation}
h_\ell\,s_\nu^{(k)}
= \sum_{\mu} c_{\mu\nu}\, s_\mu^{(k)}
\,,
\end{equation}
for some coefficients $c_{\mu\nu}$.  To determine the
$c_{\mu\nu}$, we examine
$h_\ell h_\lambda$.  Using the $k$-Schur expansion \eqref{e1} 
for $h_\lambda$, we find that
\begin{equation}
\label{pierieq2}
h_{\ell} h_\lambda=
\sum_{\nu } K^{(k)}_{\nu\lambda} \,h_\ell\,s_\nu^{(k)}
=
\sum_{\nu } K^{(k)}_{\nu\lambda} 
\sum_{\mu} c_{\mu\nu}\, s_\mu^{(k)}
\,.
\end{equation}
On the other hand, we can use \eqref{e1} to expand $h_{\ell} h_\lambda=
h_{(\ell,\lambda)}$.
Then applying Corollaries~\ref{corosym} and \ref{kostka2} we obtain
\begin{equation}
\label{pierieq1}
h_{(\ell,  \lambda)}=
\sum_{\mu} K^{(k)}_{\mu\,(\ell,\lambda)} \,s_\mu^{(k)} 
=
\sum_{\mu } 
\sum_{{\mu/\nu=\text{horizontal $\ell$-strip}}\atop
\mu^{\omega_k}/\nu^{\omega_k}=\text{vertical $\ell$-strip}}
\!\!\! \!\!\!  \!\!\! \!\!\!
K^{(k)}_{\nu\lambda} \,s_\mu^{(k)}
\,.
\end{equation}
We can equate the coefficient of $s_\mu^{(k)}$ in the right side of this 
expression to that of \eqref{pierieq2} to get the system:
\begin{equation}
\sum_{{\mu/\nu=\text{horizontal $\ell$-strip}}\atop
\mu^{\omega_k}/\nu^{\omega_k}=\text{vertical $\ell$-strip}}
\!\!\! \!\!\!
\!\!\! \!\!\!
K^{(k)}_{\nu\lambda} 
= \sum_{\nu} K^{(k)}_{\nu\lambda} c_{\mu\nu}
\,.
\end{equation}
One obvious solution is the one we want:
$$
c_{\mu \nu} = 
\left\{ 
\begin{array}{ll} 
1 & {\rm if~} \mu \in H_{\nu,\ell}^{(k)} \\
0 & {\rm otherwise}
\end{array}
\right\}
\,.
$$
In fact, this is the unique solution since another solution 
$c_{\mu \nu}'$ would satisfy
\begin{equation}
0= \sum_{\nu} K^{(k)}_{\nu\lambda} ( c_{\mu \nu}'-c_{\mu\nu})
\, ,
\end{equation}
where the invertibility of the matrix
$K^{(k)}_{\nu\lambda}$ implies $c_{\mu \nu}'=c_{\mu\nu}$.
\end{proof}

\medskip

Skew $k$-tableaux can be used to encode the iteration 
of the $k$-Pieri rule, generalizing Theorem~\ref{kpieri}.

\medskip

\begin{definition}
\label{defktabgenskew} 
Let $\delta\subseteq\gg$ be $k+1$-cores with $m_1$ and $m_2$
$k$-bounded hooks respectively, and  let
$\aa=(\aa_1,\ldots,\aa_r)$ be a composition of $m_1-m_2$.  A 
\dit{skew $k$-tableau} of shape $\gg/\delta$ and ``$k$-weight" $\aa$ 
is a filling of $\gg/\delta$ with integers $1,2,\ldots,r$ such that

\smallskip
\noindent
(i) rows are weakly increasing and columns are strictly increasing

\smallskip
\noindent
(ii) the collection of cells filled with letter $i$ are labeled by exactly 
$\alpha_i$ distinct $k+1$-residues.
\end{definition}

\begin{remark}
\label{extskew}
Our results on $k$-tableaux easily extend to include skew $k$-tableaux.
In particular, the discussion in \S\ref{ktabprops} of how to obtain a
$k$-tableau by deleting the largest letters from a given $k$-tableau 
explains more generally that deleting the largest letter from
a skew $k$-tableau produces a valid skew $k$-tableau.  Furthermore, 
although we have defined $\tau_a$ on $k$-tableaux, the results clearly
hold for skew $k$-tableaux.
\end{remark}

\begin{corollary} 
\label{propcons1}
For any $k$-bounded partitions $\lambda$ and $\mu$,
\begin{equation} 
\label{hinduhyp}
h_{\lambda}\, s_{\mu}^{(k)} = \sum_{\nu\in\mathcal P^k} 
K^{(k)}_{\nu/\mu,\lambda} \, s_{\nu}^{(k)} \, ,
\end{equation}
where $K^{(k)}_{\nu/\mu,\lambda}$ is the number of skew $k$-tableaux of 
shape $\core(\nu)/\core(\mu)$ and $k$-weight $\lambda$.
\end{corollary}
\begin{proof}  
If $\nu\in H_{\mu,\ell}^k$ then $\core(\nu)/\core(\mu)$ is a
horizontal strip with $\ell$ residues by \eqref{remark23}.  
Thus, the $k$-Pieri rule implies our claim when $\lambda=(\ell)$ 
and we proceed by induction on $\ell(\lambda)$.
Assuming \eqref{hinduhyp} holds for $\lambda$ with $\ell(\lambda)<n$,
we have
\begin{equation}
h_{(\ell,\lambda)}\, s_{\mu}^{(k)} 
=h_{\ell} h_{\lambda}\, s_{\mu}^{(k)} 
= \sum_{\nu} K^{(k)}_{\nu/\mu,\lambda} \, h_{\ell} \, s_{\nu}^{(k)} 
= \sum_{\omega} \sum_{\nu} K^{(k)}_{\nu/\mu  ,\, \lambda} \,   
K^{(k)}_{\omega/\nu , (\ell)  } \, s_{\omega}^{(k)} \, .
\end{equation}
Since removing the highest letter from a skew $k$-tableau produces a 
skew $k$-tableau by Remark~\ref{extskew}, we have that
\begin{equation}
\sum_{\nu} K^{(k)}_{\nu/\mu, \,  \lambda} \,   K^{(k)}_{\omega/\nu , (\ell)  } = 
K^{(k)}_{\omega/\mu  , \, (\ell,\lambda)} \, ,
\end{equation}
implying our claim.
\end{proof}

\medskip

As with the Schur functions, there is also a combinatorial rule to 
compute $e_\ell s_\nu^{(k)}$ in terms of $k$-Schurs by
using vertical strips to $\nu$ rather than horizontal.

\begin{theorem}
\label{ekpieri}
For any $k$-bounded partition $\nu$ and integer $\ell\leq k$, 
\begin{equation}
\label{epieri}
e_\ell\,s_\nu^{(k)} = \sum_{ \lambda\in E_{\nu,\ell}^{(k)}}
s_\lambda^{(k)} \,,
\end{equation}
where $E_{\nu,\ell}^{(k)} =
\left\{ \lambda \, \Big| \, 
\lambda/\nu={\rm vertical~} \ell\text{-}{\rm strip} \;\; {\rm and}
\;\;
\lambda^{\omega_k}/\nu^{\omega_k}={\rm horizontal~} 
\ell\text{-}{\rm strip} \right\} \,.$
\end{theorem}

\medskip

In this case, $\lambda\in E_{\nu,\ell}$ implies 
$\core(\lambda^{\omega_k})/\core(\nu^{\omega_k})$ is a vertical
strip with $\ell$ distinct residues by \eqref{remark23}.  We can
thus apply the same argument used to derive Corollary~\ref{propcons1}
from Theorem~\ref{kpieri} to prove the corollary:

\medskip

\begin{corollary} \label{corocons2}
For any $k$-bounded partitions $\lambda$ and $\mu$,
\begin{equation} 
\label{induhyp}
e_{\lambda}\, s_{\mu}^{(k)} = \sum_{\nu} 
\tilde K_{\nu/\mu,\lambda} \, s_{\nu}^{(k)} \, ,
\end{equation}
where $\tilde K_{\nu/\mu,\lambda}$ is the number of 
\dit{transposed skew $k$-tableaux} of shape $\core(\nu)/\core(\mu)$ 
and $k$-weight $\lambda$.  Such tableaux are defined by the same 
conditions as skew $k$-tableaux except that condition (i) is 
changed to: rows are {strictly} increasing and columns are 
{weakly} increasing.
\end{corollary}

\medskip

\noindent
{\it Proof of Theorem~\ref{ekpieri}.}
Since $e_1=h_1$, Theorem~\ref{kpieri} implies the case when $\ell=1$ and
we assume by induction that the action of $e_r$ for all $r<\ell$ is given 
by \eqref{epieri}.   To prove our assertion for multiplication by
$e_\ell$, note that by applying the identity (e.g. \cite{[Ma]}):
$$
\sum_{r=0}^{\ell-1} (-1)^r h_{\ell-r} \,e_r  + (-1)^{\ell}\, e_\ell = 0 \,,
$$
Eq.~\eqref{epieri} follows from the expression
\begin{equation} 
\label{epieri2}
\sum_{r=0}^{\ell-1} (-1)^r h_{\ell-r} \,e_r \,s_\nu^{(k)}
+
(-1)^{\ell}
\sum_{\lambda\in E_{\nu,\ell}^{(k)}} s_\lambda^{(k)} = 0 \,.
\end{equation}
It thus suffices to show the coefficient of $s_\mu^{(k)}$ in the 
left side of this expression is zero.

By induction, Corollaries~\ref{propcons1} and \ref{corocons2} 
tell us that for $r<\ell$, the coefficient of
$s_{\mu}^{(k)}$ in $h_{\ell-r}e_r s_\nu^{(k)}$
is the number of fillings with weight $(r,\ell-r)$ in the
following set:

\begin{definition}
Let $\mu,\nu$ be $k$-bounded partitions and fix $\ell\leq k$.  An 
element $T\in A^{(k)}_{\nu,\ell}(\mu)$ of weight $(r,\ell-r)$, 
for any $0\leq r\leq \ell$, has shape $\core(\mu)/\core(\nu)$
and is filled with letters $x<y$ such that
\begin{enumerate}
\item[(i)]  $T_{\leq x}$ is a {transposed} skew $k$-tableau
of $k$-weight $(r)$ filled with letter $x$

\item[(ii)] $T/T_{\leq x}$ is a skew $k$-tableau
of $k$-weight $(\ell-r)$ filled with letter $y$
\end{enumerate}
\end{definition}

Since the the coefficient of $s_{\mu}^{(k)}$ in 
$\sum_{\lambda\in E_{\nu,\ell}^{(k)}} s_\lambda^{(k)}$
is the number of fillings in $\mathcal A_{\nu,\ell}^{(k)}(\mu)$
with weight $(\ell,0)$, the coefficient of 
$s_{\mu}^{(k)}$ in the left side of \eqref{epieri2} equals
\begin{equation*}
\sum_{r=0}^{\ell-1}
\!\! \!\!
\sum_{T\in \mathcal A^{k}_{\nu,\ell}(\mu)\atop
weight(T)=(r,\ell-r)}
\!\! \!\!  \!\! \!\!
(-1)^{r} \quad + \sum_{T\in \mathcal A^{k}_{\nu,\ell}(\mu)\atop
weight(T)=(\ell,0)}
\!\! \!\!  \!\! \!\!
 (-1)^\ell \quad = \quad  
\sum_{T\in \mathcal A^{k}_{\nu,\ell}(\mu)}
\!\! \!\!
(-1)^{{\rm sgn}(T)} \, ,
\end{equation*}
where ${\rm sgn}(T)=(-1)^r$ for weight$(T)=(r,\ell-r)$.
If we can produce a sign reversing involution 
$\mathfrak m$ on $\mathcal A^{k}_{\nu,\ell}(\mu)$, then
$$
\sum_{T\in \mathcal A^{k}_{\nu,\ell}(\mu)}
\!\! \!\!
(-1)^{{\rm sgn}(T)} = \sum_{T\in \mathcal A^{k}_{\nu,\ell}(\mu)}
\!\! \!\!
(-1)^{{\rm sgn}\bigl(\mathfrak m(T) \bigr)}= 
-
\!\! \!\!
\sum_{T\in \mathcal A^{k}_{\nu,\ell}(\mu)}
\!\! \!\!
(-1)^{{\rm sgn}(T)} \, ,
$$
implying the coefficient of $s_{\mu}^{(k)}$ is zero.
Proposition~\ref{minvo} below gives the desired $\mathfrak m$.
\hfill$\square$

\medskip

The involution $\mathfrak m$ acts on \dit{free} entries
of $T\in \mathcal A_{\nu,\ell}^{(k)}(\mu)$, where
an $x(i)$ is free if every $x(i)\in T$ occurs at 
the top of its column, and $y(j)$ is free if no
$y(j)$ is right-adj to an $x$ or $y$.

\begin{definition}
The map $\mathfrak m$ acts on $T\in\mathcal A_{\nu,\ell}^{(k)}(\mu)$ by: 

\noindent
1) Let $r_1$ denote the lowest row containing a free $x$ and $i$ 
denote its residue (if there is no free $x$, set $r_1=\infty$).
Let $r_2$ be the lowest row containing a free $y$ and $j$ its residue
(if there is no free $y$, set $r_2=\infty$).

\smallskip

\noindent
2)  If $r_1<r_2$, send every $x(i)$ to $y(i)$.  Otherwise
send every $y(j)$ to $x(j)$. 
\end{definition}

The definition of $\mathfrak m$ is well-defined since
every $T\in \mathcal A_{\nu,\ell}^{(k)}(\mu)$ contains 
a free $x$ or a free $y$.  For example, $x(i)$ is not free 
in $T$ implies there is an $x(i-1)$ or a $y(i-1)$ in $T$ and
$y(i)$ is not free implies there is an $x(i-1)$ or a $y(i-1)$ 
in $T$.  By iteration, no letter is free implies that $T$ 
contains
$$
z(i) , z(i-1), z(i-2),\dots,z(i+2),z(i+1)\,,
$$
with each $z(j)=x(j)$ or $y(j)$.  This contradicts that
$T$ has weight $(r,\ell-r)$ for $\ell\leq k$.

\begin{proposition}
\label{minvo}
The map $\mathfrak m$ is an involution on $\mathcal A_{\nu,\ell}^{(k)}(\mu)$ 
such that weight$({\mathfrak m(T)})=({n_1\pm 1},n_2\mp 1)$, given
weight$(T)=({n_1},n_2)$.  
\end{proposition}
\begin{proof}
Let $\hat T=\mathfrak m(T)$ for some 
$T\in\mathcal A_{\nu,\ell}^{(k)}(\mu)$.  First note that the 
definition of free implies the $x$'s form a vertical strip 
and the $y$'s a horizontal strip in $\hat T$.  

To determine
how the weight of $T$ changes under $\mathfrak m$, consider
first the case that $r_1<r_2$.  
Since every $x(i)\in T$
goes to $y(i)\in \hat T$, there are only $n_1-1$ residues 
of $x$ in $\hat T$.  To show that there are $n_2+1$ residues
of $y$ in $\hat T$, we must prove $y(i)\not\in T$.
Suppose there is a $y(i)$ in $T$.  Since $T/T_{\leq x}$ is a skew 
$k$-tableau, $T_{<y(i)}/T_{\leq x}$ is a skew $k$-tableaux by 
Remark~\ref{extskew}.  Thus, $T_{<y(i)}$ has core shape
and an addable corner $y(i)$ of residue $i$.
Further, $x(i)$ is a removable corner in $T_{<y(i)}$ since
$x(i)$ is at the top of its column and $y(i+1)\not\in T_{<y(i)}$.
We reach a contradiction by Remark~\ref{noboth} which tells
us a core cannot have an addable and removable corner of
the same residue.  A similar argument works when $r_2<r_1$.

Lastly, to see that $\mathfrak m$ is an involution, consider the case 
that $r_1<r_2$.  Since there is at most one $x$ in each row, any row 
where $\mathfrak m: x(i)\to y(i)$ contains only $y$'s in $\hat T$.  
Thus $y(i)$ is free in $\hat T$ and $r_1$ is the lowest row with $y(i)$
since $y(i)\not\in T$ by the previous paragraph.  There are
no lower free entries by definition of $r_1$.  Therefore, when
$\mathfrak m$ is applied to $\hat T$, $y(i)\to x(i)$ and $T$ is 
recovered.  Similarly when $r_1>r_2$.
\end{proof}

\subsection{Further properties}

Recall the algebra endomorphism $\omega$ that provides an involution
on $\Lambda$,  defined by $\omega h_\ell = e_\ell$.  This map has an 
especially simple action on the Schur functions: $\omega s_\lambda=s_{\lambda'}$.  
Since $\omega$ is also an involution on $\Lambda^{(k)}$, we can ask how it
acts on a $k$-Schur function.  Right on cue, we find:

\medskip

\begin{theorem}
The $\omega$-involution acts on the $k$-Schur functions by
\begin{equation}
\omega s_\lambda^{(k)} = s_{\lambda^{\omega_k}}^{(k)} \, .
\end{equation}
\end{theorem}

\begin{proof}
Let $F_\mu= \omega s_{\mu^{\omega_k}}^{(k)}$.  Since
$h_\ell \, \omega\left(s_{\lambda^{\omega_k}}^{(k)}\right)
=\omega\left(e_\ell s_{\lambda^{\omega_k}}^{(k)}\right)$, 
we can apply the $k$-Pieri rule (Theorem~\ref{ekpieri}) to obtain
\begin{equation}
h_\ell\,F_\lambda=
\omega\left(e_\ell \, s_{\lambda^{\omega_k}}^{(k)}\right)
= \sum_{\mu \in E_{\lambda^{\omega_k},\ell}^{(k)} }
\omega\, s_\mu^{(k)} 
\;=\; 
\sum_{\mu^{\omega_k} \in E_{\lambda^{\omega_k},\ell}^{(k)} }
F_{\mu}
\;=\; \sum_{\mu \in H_{\lambda,\ell}^{(k)}}
F_{\mu}
\, ,
\end{equation} 
recalling that $(\mu^{\omega_k})^{\omega_k}=\mu$.
Iteration of this expression from $F_0=\omega s_0^{(k)}=1$ matches
iteration of the $k$-Pieri rule from $s_0^{(k)}=1$.  Thus, $F_\mu$ 
satisfies 
\begin{equation}
h_\lambda = F_\lambda + \sum_{\mu:\mu\rhd\lambda}
K_{\mu\lambda}^{(k)}\, F_\mu
\end{equation}
implying that $F_\mu=s_\mu^{(k)}$ by Definition~\ref{kschurdef}
of the $k$-Schur functions.
\end{proof}

\medskip

 From the action of $\omega$ on a $k$-Schur function, we are able to
show that a $k$-Schur function reduces simply to a Schur function
when $k$ is large.  
\begin{property}
\label{kss}
For any partition $\lambda$ with main hook-length  $h(\lambda)\leq k$,
we have that $s_\lambda^{(k)}=s_\lambda$.
\end{property}
\begin{proof}
Given the triangular form \eqref{triom},
\begin{equation}
\label{lowtrip}
s_{\lambda}^{(k)}=s_{\lambda} +higher~terms\,,
\end{equation}
we can apply the $\omega$-involution to obtain:
\begin{equation}
\label{lowtri}
s_{\lambda^{\omega_k}}^{(k)}=s_{\lambda'} +lower~terms\,.
\end{equation}
However, since $\lambda^{\omega_k}=\lambda'$ when $h(\lambda)\leq k$
from Remark~\ref{conjupetit}, 
the previous expression reduces to
\begin{equation}
s_{\lambda'}^{(k)}=s_{\lambda'} +lower~terms\,.
\end{equation}
Setting this equal to \eqref{lowtrip}, with
$\lambda$ replaced by $\lambda'$,
proves our claim.
\end{proof}

We finish by deriving one last property from the action of the 
$\omega$-involution.  This property is one of the few that we were
able to prove using a prior characterization (see \cite{[LMadv]}).  
In particular, there exists a subset of \dit{irreducible} $k$-Schur functions 
from which all other $s_\lambda^{(k)}$ may be constructed with multiplication 
by usual Schur functions indexed by \dit{$k$-rectangles} -- partitions of the 
form $(\ell^{k-\ell+1})$.  The irreducibles 
consist of the special set of $k$-Schur functions indexed by irreducible 
partitions; $k$-bounded partitions with no more than $i$ parts equal 
to $k-i$, for $i=0,\dots,k-1$.  Remarkably, we can also prove this result
using the characterization studied in this article.

\begin{theorem}
For any $k$-rectangle $\square$ and $k$-bounded partition $\mu$, we have
$$
s_\square s_\mu^{(k)}= s_{\mu \cup \square}^{(k)}
\,.
$$
\end{theorem}
\begin{proof}
Consider the linear operator $\Theta_\square$ defined on $\Lambda^{(k)}$ 
by $\Theta_\square s_\mu^{(k)} = s_{\mu\cup\square}^{(k)}$.
It suffices to show that 
$\Theta_\square s_\mu^{(k)}=s_\mu^{(k)}\Theta_\square$ since
$\Theta_\square\cdot 1= \Theta_\square s^{(k)}_\emptyset = 
s^{(k)}_\square=s_\square$ by Property~\ref{kss}.  However, 
since the homogeneous functions generate $\Lambda^{(k)}$, we
will instead prove that
$\Theta_\square h_\ell = h_\ell \Theta_\square$.
To this end, note that the $k$-Pieri rule \eqref{kpierieq} implies 
\begin{eqnarray} 
\label{eq37}
\Theta_\square  h_\ell \, s_\mu^{(k)} 
& =  & \Theta_\square 
\sum_{ \eta \in H_{\mu,\ell}^{(k)}}
s_\eta^{(k)}
= 
\sum_{ \eta \in H_{\mu,\ell}^{(k)}}
s_{\eta\cup\square}^{(k)}\, ,
\end{eqnarray}
and on the other hand,
\begin{eqnarray}
\label{intrec}
h_\ell \, \Theta_\square s_\mu^{(k)} & = & h_\ell \, 
s_{\mu\cup\square}^{(k)} = 
\sum_{ \gamma \in H_{\mu \cup \square,\ell}^{(k)}}
s_{\gamma}^{(k)} \, .
\end{eqnarray}
It is known (Corollary 57 in \cite{[LMcore]}) that 
$\gamma \in H_{\mu \cup \square,\ell}^{(k)}$
implies $\mu\cup\square\preceq\gamma$, where
$\alpha\preceq\beta$ is defined on $k$-bounded partitions
by the covering relation: $\alpha\prec\!\!\!\cdot\,\beta$ when 
$\beta,\alpha$ are 1-admissible.
Then, using Theorem 20 from \cite{[LMrec]}:
$\mu\cup\square\preceq\gamma\iff\gamma=\eta\cup\square\,
\text{ and }\, \mu\preceq \eta$ for some $k$-bounded $\eta$, we 
can transform \eqref{intrec} into
\begin{equation} \label{eq39}
h_\ell\,\Theta_\square s_\mu^{(k)} =  
\sum_{ \eta \cup \square \in H_{\mu \cup \square,\ell}^{(k)}}
s_{\eta \cup \square}^{(k)} \, .
\end{equation}

Since the $k$-Schur functions form a basis for $\Lambda^{(k)}$,
it remains to show that the right side of
Eq.~\eqref{eq37} equals that of Eq.~\eqref{eq39}, or equivalently
that
$$
\eta \cup \square \in H_{\mu \cup \square,\ell}^{(k)}\iff
\eta  \in H_{\mu,\ell}^{(k)}\,.
$$
Given $(\eta\cup\square)^{\omega_k}=\eta^{\omega_k}\cup\square^{\omega_k}$
by Theorem~10 of \cite{[LMrec]}, we have
$\eta \cup \square \in H_{\mu \cup \square,\ell}^{(k)}$ iff
$\eta \cup \square/ \mu \cup \square$ is a horizontal strip and
$\eta^{\omega_k} \cup \square'/\mu^{\omega_k}\cup\square'$ is 
vertical strip.  Thus, our claim follows by noting that
$\alpha \cup \square/ \beta \cup \square$ is a horizontal (resp.
vertical) strip iff $\alpha/\beta$ is a horizontal (resp. vertical) strip.
\end{proof}

\bigskip

\noindent{\bf Acknowledgments}
{\it We thank Michelle Wachs for her ideas and her suggestion to explore the 
fruitful characterization used here.}

\end{document}